\documentclass[10pt]{article}

\usepackage[dvips,letterpaper,top=0.5in, bottom=0.5in, left=0.75in, right=0.5in,includefoot,heightrounded]{geometry}

\usepackage{srcltx}

\usepackage{amsmath}
\usepackage{amssymb,amsfonts}

\usepackage{abstract}

\usepackage{epsfig}
\usepackage[usenames,dvipsnames]{color}
\usepackage[usenames,dvipsnames]{xcolor}

\usepackage{booktabs}

\usepackage{flushend}
\usepackage{url}

\usepackage{enumerate}

\usepackage[short,12hr]{datetime}

\usepackage{hyperref}

\usepackage{lipsum}

\usepackage[sc,tiny]{titlesec}
       \usepackage{mathptmx}

\newtheorem{definition}{Definition}
\newtheorem{theorem}{Theorem}



\title{%
Fast Finite Field Hartley Transforms Based on Hadamard Decomposition}

\author{%
H. M. de Oliveira%
\thanks{%
H. M. de Oliveira
was with 
Departamento de Eletr\^onica e Sistemas,
Universidade Federal de Pernambuco (UFPE).
He is currently
with
the Signal Processing Group,
Departamento de Estat\'{\i}stica, 
Universidade Federal de Pernambuco.
Email: \url{hmo@ufpe.br}
}
\and
R. G. F. T\'avora%
\thanks{%
R. G. F. T\'avora
was with 
the Communications Research Group,
Departamento de Eletr\^onica e Sistemas,
Universidade Federal de Pernambuco.
Email: \url{r_tavora@hotmail.com}
}
\and
R. J. Cintra%
\thanks{%
R. J. Cintra 
was with 
the Communications Research Group,
Departamento de Eletr\^onica e Sistemas,
Universidade Federal de Pernambuco.
He is currently
with
the Signal Processing Group,
Departamento de Estat\'{\i}stica, 
Universidade Federal de Pernambuco.
E-mail: \protect\url{rjdsc@de.ufpe.br},
}
\and
R. M. Campello de Souza%
\thanks{%
R. M. Campello de Souza
is with 
the Communications Research Group,
Departamento de Eletr\^onica e Sistemas,
Universidade Federal de Pernambuco.
Email: \url{ricardo@ufpe.br}
}
}

\date{\today\ @ \currenttime}
\date{}

\newcommand{\myabstract}{%
A new transform over finite fields, 
the 
finite field Hartley transform (FFHT), 
was recently introduced and 
a number of promising applications on the design of 
efficient multiple 
access systems and multilevel spread spectrum
sequences were proposed.
The FFHT exhibits interesting symmetries, which 
are exploited to derive tailored fast transform algorithms.
The proposed fast algorithms are
based on successive decompositions of the FFHT
by means of Hadamard-Walsh transforms (HWT).
The introduced decompositions
meet the lower bound on the multiplicative complexity 
for all the cases investigated.
The complexity of the new algorithms
is compared with that of traditional algorithms.
}

\newcommand{\mykeywords}{%
Finite field transforms,
fast algorithms,
discrete Hartley transform
}

\begin{document}

\twocolumn[%
  \maketitle
  \begin{onecolabstract}
    \myabstract
  \end{onecolabstract}
  \begin{center}
    \small
    \textbf{Keywords~\vspace{3mm}}
    \linebreak
    \mykeywords
  \end{center}
  \bigskip
]
\saythanks

\section{Introduction}

Discrete transforms defined over finite fields,
such as 
the finite field Fourier transform (FFFT),
pivotal tools in coding theory~\cite{BLAH} and signal processing~\cite{Reed}.
Another interesting example is the 
finite field Hartley transform
(FFHT), 
a self-inverse transform  (involution operator)
introduced in~\cite{TRIGO,FFHT,CFFHT}.
Recent promising applications of discrete
transforms concern the use of the FFHT to design digital multiplex systems,
efficient multiple access systems~\cite{MA} and multilevel spread spectrum
sequences~\cite{Spread}.
A decisive factor for applications of discrete
transforms has been the existence of the so-called fast transforms (FT) for
computing it.
Since the FFHT is a more symmetrical version of discrete
transform, in this paper this symmetry is exploited so as to derive new FTs that require
less operations.
These FTs, derived for short blocklengths ($N\leq 24)$, are
based on successive decompositions in a similar way as the multilayer
Hadamard decomposition employed~\cite{Renato} 
to compute
the discrete Hartley Transform (DHT)~\cite{BRAC}.
This new approach,
which is based on
decomposition of the FFHT by means of Hadamard-Walsh transforms (HWT), 
meets the lower bound on the multiplicative complexity of
a discrete Fourier transform (DFT)~\cite{Complexid}.
Each
HWT implements pre-additions and post-additions.
These schemes are easy to
implement using digital signal processors (DSP) or low-cost high-speed
dedicated hardware.
The complexity of these new FTs is compared
with that of traditional methods, 
such as the 
Cooley-Tukey radix-2,
split radix, 
Winograd, 
and 
Rader-Brenner algorithms, 
which were adapted to compute
the FFHT~\cite{Chile}. 

\section{The Finite Field Hartley Transform}

Finite field Hartley transforms are
based on a trigonometry over Galois Fields $GF(q)$, $q=p^r$, $p\equiv3\pmod4$,
so that $(p-1)^{1/2}\not \in GF(q)$. The set $G(q)$ of
Gaussian integers over $GF(q)$\ plays an important role in this analysis.
This set defines a structure $GI(q)$,
which is  isomorphic to 
$GF(q^2)$~\cite{TRIGO}.

\begin{definition}
Let $\zeta $ be an element of $GI(q)$ with multiplicative order $N$, where 
$q=p^r$. The trigonometric functions sine, cosine, and 
$\operatorname{cas}$ (cosine-and-sine or Hartley kernel)
are defined by,
respectively: 
\begin{align*}
\sin (i) &= \frac{\zeta ^i-\zeta ^{-i}}{2j},
\\
\cos (i)&=\frac{\zeta ^i+\zeta ^{-i}}2,
\\
\operatorname{cas}(i)&=\sin (i)+\cos (i)
,
\end{align*}
for $i = 0,1,\ldots,N-1$.
\end{definition}

\begin{definition}
Let $v=\{v_0,v_1,\ldots,v_{N-1}\}$ be a vector of $GF(q)$-valued components, 
$q=p^r$.
The finite field Hartley transform (FFHT) is the vector 
$V=\{V_0,V_1,\ldots,V_{N-1}\}$, 
with components $V_k\in GI(q^m)$ given by
$V_k=\sum_{i=0}^{N-1}v_i\cdot\operatorname{cas}(ik)$,
where $\zeta $ is an element of multiplicative order $N$ over $GI(q^m)$.
\end{definition}
The inverse FFHT is given by the following theorem.
\begin{theorem}
The vector $v=\{v_0,v_1,\ldots,v_{N-1}\}$ can be derived from its FFHT according
to: 
\begin{align*}
v_i
=
\frac{1}{N\pmod{p}}
\sum_{k=0}^{N-1}
V_k\cdot \operatorname{cas}(ik)
,
\end{align*}
for $i = 0,1,\ldots,N-1$.
\end{theorem}

\section{Hadamard Decomposition of the FFHT}

The Hadamard decomposition was employed 
in~\cite{Renato} as a tool to compute the discrete Hartley transform.
This approach
allows the minimization of the multiplicative complexity of the DHT for some
blocklengths.
Since all the properties and symmetries of the DHT are also
observed for the FFHT, the application of this algorithm to finite fields
should be expected.
The minimal multiplicative complexity
of a DFT
with blocklength $N$---denoted by $\mu (DFT(N))$---can be calculated by converting the DFT in a set of
cyclic convolutions. A lower bound on $\mu (DFT(N))$ is presented in~\cite{Complexid}. 
Table~\ref{mi} shows a few values of $\mu (DFT(N))$ for short
blocklengths.

\begin{table}
\begin{center}
\caption{Minimal Multiplicative Complexity achievable for the $N$-point DFT}
\begin{tabular}{ccccc}
\toprule 
$N$ & 4 & 8 &12 & 16  \\
\midrule
$\mu(DFT(N))$  & 0 & 2 & 4 & 10 \\
\bottomrule
\end{tabular}

\label{mi}
\end{center}
\end{table}

Considering with finite field transforms, 
the following comments are worthwhile:

\begin{enumerate}[(i)]

\item  The minimal multiplicative complexity, $\mu (FFFT(N))$, for a FT over
the finite field $GI(p^r)$, is the same as $\mu (DFT(N)),$ evaluated over
the real field.

\item  The relationship between the multiplicative and additive complexity
over a finite field strong depends on implementation. For small $p,$ the
total complexity (additive plus multiplicative) must be taken into account
since their difference is small.
\end{enumerate}
New algorithms for computing the FFHT are introduced in the next section.

\subsection{Computing the 4-point FFHT}

Let $v\longleftrightarrow V$ be a FFHT transform pair over $GI(7)$. The
FFHT, assuming a $\operatorname{cas}(\cdot)$ kernel with $\zeta =j$, is computed by: 
$$
\left[ 
\begin{smallmatrix}   
V_0 \\ 
V_1 \\ 
V_2 \\ 
V_3 
\end{smallmatrix}
\right] =\left[ 
\begin{smallmatrix}
1 & 1 & 1 & 1 \\ 
1 & 1 & 6 & 6 \\ 
1 & 6 & 1 & 6 \\ 
1 & 6 & 6 & 1 
\end{smallmatrix}
\right] \left[ 
\begin{smallmatrix}   
v_0 \\ 
v_1 \\ 
v_2 \\ 
v_3 
\end{smallmatrix}
\right] . 
$$
Indeed, no multiplication is needed. Observing further symmetries, columns
can be combined through Hadamard blocks in order to reduce the number of
additions. Let 
\begin{align*}
S_0(1)&=(v_3-v_1),\,S_1(1)=(v_3+v_1),\\
S_2(1)&=(v_0-v_2),\,S_3(1)=(v_0+v_2).
\end{align*}
It follows that:
$$
\left[ 
\begin{smallmatrix}
V_0 \\ 
V_1 \\ 
V_2 \\ 
V_3 
\end{smallmatrix}
\right] =\left[ 
\begin{smallmatrix}
0 & 1 & 0 & 1 \\ 
6 & 0 & 1 & 0 \\ 
0 & 6 & 0 & 1 \\ 
1 & 0 & 1 & 0 
\end{smallmatrix}
\right] \left[ 
\begin{smallmatrix}
S_0(1) \\ 
S_1(1) \\ 
S_2(1) \\ 
S_3(1) 
\end{smallmatrix}
\right] . 
$$
The number of additions is reduced from 12 to 8 (4 pre-additions and 4
post-additions).

\subsection{Computing the $6$-point FFHT}

Let $v\longleftrightarrow V$ be an FFHT transform pair over $GI(7)$.
Considering $\zeta =3$, the FFHT can be computed by

$$
\left[ 
\begin{smallmatrix}
V_0 \\ 
V_1 \\ 
V_2 \\ 
V_3 \\ 
V_4 \\ 
V_5 
\end{smallmatrix}
\right] =\left[ 
\begin{smallmatrix}
1 & 1 & 1 & 1 & 1 & 1  \\ 
1 & 4+j & 3+j & 6 & 3+6j & 4+6j \\ 
1 & 3+j & 3+6j & 1 & 3+j & 3+6j \\ 
1 & 6 & 1 & 6 & 1 & 6 \\ 
1 & 3+6j & 3+j & 1 & 3+6j & 3+j \\ 
1 & 4+6j & 3+6j & 6 & 3+j & 4+j 
\end{smallmatrix}
\right] \left[ 
\begin{smallmatrix}
v_0 \\ 
v_1 \\ 
v_2 \\ 
v_3 \\ 
v_4 \\ 
v_5 
\end{smallmatrix}
\right] . 
$$
Observing the symmetries, a first column combination can be made. Let 
$$
\begin{smallmatrix}
S_0(1)=(v_4-v_1),\,S_1(1)=(v_4+v_1),\,S_2(1)=(v_5-v_2), \\ 
S_3(1)=(v_5+v_2),\,S_4(1)=(v_0-v_3),\,S_5(1)=(v_0+v_3). 
\end{smallmatrix}
$$
Therefore, 
$$
\left[ 
\begin{smallmatrix}
V_0 \\ 
V_1 \\ 
V_2 \\ 
V_3 \\ 
V_4 \\ 
V_5 
\end{smallmatrix}
\right] =\left[ 
\begin{smallmatrix}
0 & 1 & 0 & 1 & 0 & 1 \\ 
3+6j & 0 & 4+6j & 0 & 1 & 0 \\ 
0 & 3+j & 0 & 3+6j & 0 & 1 \\ 
1 & 0 & 6 & 0 & 1 & 0 \\ 
0 & 3+6j & 0 & 3+j & 0 & 1 \\ 
3+j & 0 & 4+j & 0 & 1 & 0 
\end{smallmatrix}
\right] \left[ 
\begin{smallmatrix}
S_0(1) \\ 
S_1(1) \\ 
S_2(1) \\ 
S_3(1) \\ 
S_4(1) \\ 
S_5(1) 
\end{smallmatrix}
\right] . 
$$
Going on with this procedure, a second pre-addition layer is derived: 
\begin{align*}
S_0(2)&=S_2(1)-S_0(1),  S_1(2)=S_2(1)+S_0(1), \\
S_2(2)&=S_3(1)-S_1(1), S_3(2)=S_3(1)+S_1(1), \\
S_4(2)&=S_4(1), S_5(2)=S_5(1)
.
\end{align*}
Finally, 
$$
\left[ 
\begin{smallmatrix}
V_0 \\ 
V_1 \\ 
V_2 \\ 
V_3 \\ 
V_4 \\ 
V_5 
\end{smallmatrix}
\right] =\left[ 
\begin{smallmatrix}
0 & 0 & 0 & 1 & 0 & 1 \\ 
4 & 6j & 0 & 0 & 1 & 0 \\ 
0 & 0 & 6j & 3 & 0 & 1 \\ 
6 & 0 & 0 & 0 & 1 & 0 \\ 
0 & 0 & j & 3 & 0 & 1 \\ 
4 & 1j & 0 & 0 & 1 & 0 
\end{smallmatrix}
\right] \left[ 
\begin{smallmatrix}
S_0(2) \\ 
S_1(2) \\ 
S_2(2) \\ 
S_3(2) \\ 
S_4(2) \\ 
S_5(2) 
\end{smallmatrix}
\right] . 
$$
Since there is only one multiplication (by the same factor) in columns 1 and
4, there will be two multiplications. The total number of additions required
to compute a 6-blocklength FFHT is 16 (10 pre-additions and 6
post-additions).

\subsection{Computing the 8-point FFHT}

Let $v\longleftrightarrow V$ be an FFHT transform pair over $GI(7)$. 
Let 
$\zeta =2+2j$, so the corresponding matrix formulation is, 
$$
\left[ 
\begin{smallmatrix}
V_0 \\ 
V_1 \\ 
V_2 \\ 
V_3 \\ 
V_4 \\ 
V_5 \\ 
V_6 \\ 
V_7 
\end{smallmatrix}
\right] =\left[ 
\begin{smallmatrix}
1 & 1 & 1 & 1 & 1 & 1 & 1 & 1  \\ 
1 & 4 & 1 & 0 & 6 & 3 & 6 & 0  \\ 
1 & 1 & 6 & 6 & 1 & 1 & 6 & 6  \\ 
1 & 0 & 6 & 4 & 6 & 0 & 1 & 3  \\ 
1 & 6 & 1 & 6 & 1 & 6 & 1 & 6  \\ 
1 & 3 & 1 & 0 & 6 & 4 & 6 & 0  \\ 
1 & 6 & 6 & 1 & 1 & 6 & 6 & 1  \\ 
1 & 0 & 6 & 3 & 6 & 0 & 1 & 4  
\end{smallmatrix}
\right] \left[ 
\begin{smallmatrix}
v_0 \\ 
v_1 \\ 
v_2 \\ 
v_3 \\ 
v_4 \\ 
v_5 \\ 
v_6 \\ 
v_7 
\end{smallmatrix}
\right] . 
$$
Defining a 1st order pre-addition layer: 
\begin{align*}
S_0(1)&=(v_5-v_1),\,S_1(1)=(v_5+v_1), \\ 
S_2(1)&=(v_6-v_2),\,S_3(1)=(v_6+v_2), \\ 
S_4(1)&=(v_7-v_3),\,S_5(1)=(v_7+v_3), \\ 
S_6(1)&=(v_0-v_4),\,S_7(1)=(v_0+v_4). 
\end{align*}
Therefore, 
$$
\left[ 
\begin{smallmatrix}
V_0 \\ 
V_1 \\ 
V_2 \\ 
V_3 \\ 
V_4 \\ 
V_5 \\ 
V_6 \\ 
V_7 
\end{smallmatrix}
\right] =\left[ 
\begin{smallmatrix}
0 & 1 & 0 & 1 & 0 & 1 & 0 & 1 \\ 
3 & 0 & 6 & 0 & 0 & 0 & 1 & 0 \\ 
0 & 1 & 0 & 6 & 0 & 6 & 0 & 1 \\ 
0 & 0 & 1 & 0 & 3 & 0 & 1 & 0 \\ 
0 & 6 & 0 & 1 & 0 & 6 & 0 & 1 \\ 
4 & 0 & 6 & 0 & 0 & 0 & 1 & 0 \\ 
0 & 6 & 0 & 6 & 0 & 1 & 0 & 1 \\ 
0 & 0 & 1 & 0 & 4 & 0 & 1 & 0 
\end{smallmatrix}
\right] \left[ 
\begin{smallmatrix}
S_0(1) \\ 
S_1(1) \\ 
S_2(1) \\ 
S_3(1) \\ 
S_4(1) \\ 
S_5(1) \\ 
S_6(1) \\ 
S_7(1) 
\end{smallmatrix}
\right] . 
$$
Defining a 2nd pre-addition layer, 
\begin{align*}
S_0(2)&=S_0(1),\,S_1(2)=S_4(1),\\ 
S_2(2)&=S_5(1)-S_1(1),\,S_3(2)=S_5(1)+S_1(1),\\ 
S_4(2)&=S_6(1)-S_2(1),\,S_5(2)=S_6(1)+S_2(1),\\ 
S_6(2)&=S_7(1)-S_3(1),\,S_7(2)=S_7(1)+S_3(1). 
\end{align*}
Consequently, we obtain:
$$
\left[ 
\begin{smallmatrix}
V_0 \\ 
V_1 \\ 
V_2 \\ 
V_3 \\ 
V_4 \\ 
V_5 \\ 
V_6 \\ 
V_7 
\end{smallmatrix}
\right] =\left[ 
\begin{smallmatrix}
0 & 0 & 0 & 1 & 0 & 0 & 0 & 1 \\ 
3 & 0 & 0 & 0 & 1 & 0 & 0 & 0 \\ 
0 & 0 & 6 & 0 & 0 & 0 & 1 & 0 \\ 
0 & 3 & 0 & 0 & 0 & 1 & 0 & 0 \\ 
0 & 0 & 0 & 6 & 0 & 0 & 0 & 1 \\ 
4 & 0 & 0 & 0 & 1 & 0 & 0 & 0 \\ 
0 & 0 & 1 & 0 & 0 & 0 & 1 & 0 \\ 
0 & 4 & 0 & 0 & 0 & 1 & 0 & 0 
\end{smallmatrix}
\right] \left[ 
\begin{smallmatrix}
S_0(2) \\ 
S_1(2) \\ 
S_2(2) \\ 
S_3(2) \\ 
S_4(2) \\ 
S_5(2) \\ 
S_6(2) \\ 
S_7(2) 
\end{smallmatrix}
\right] . 
$$
Since there is only one multiplication (by the same factor) in columns 1 and
2, there are two multiplications. The number of additions is 22 (14
pre-additions and 8 post-additions). It is worthwhile to remark that the
additive complexity is less than the one for an 8-DFT calculation by
the Winograd algorithm.

\subsection{Computing the 12-point FFHT}

Let $v\longleftrightarrow V$ be an FFHT transform pair over $GI(p)$. As an
example, let $p=7$ and $\zeta =3j.$ Then $V=Tv$ where, 
$$
T=
\left[ 
\begin{smallmatrix}
1 & 1 & 1 & 1 & 1 & 1 & 1 & 1 & 1 & 1 & 1 & 1 \\ 
1 & 4+6j & 4+6j & 1 & 3+6j & 4+j & 6 & 3+j & 3+j & 6 & 4+j & 3+6j \\ 
1 & 4+6j & 3+6j & 6 & 3+j & 4+j & 1 & 4+6j & 3+6j & 6 & 3+j & 4+j \\ 
1 & 1 & 6 & 6 & 1 & 1 & 6 & 6 & 1 & 1 & 6 & 6 \\ 
1 & 3+6j & 3+j & 1 & 3+6j & 3+j & 1 & 3+6j & 3+j & 1 & 3+6j & 3+j \\ 
1 & 4+j & 4+j & 1 & 3+j & 4+6j & 6 & 3+6j & 3+6j & 6 & 4+6j & 3+j \\ 
1 & 6 & 1 & 6 & 1 & 6 & 1 & 6 & 1 & 6 & 1 & 6 \\ 
1 & 3+j & 4+6j & 6 & 3+6j & 3+6j & 6 & 4+6j & 3+j & 1 & 4+j & 4+j \\ 
1 & 3+j & 3+6j & 1 & 3+j & 3+6j & 1 & 3+j & 3+6j & 1 & 3+j & 3+6j \\ 
1 & 6 & 6 & 1 & 1 & 6 & 6 & 1 & 1 & 6 & 6 & 1 \\ 
1 & 4+j & 3+j & 6 & 3+6j & 4+6j & 1 & 4+j & 3+j & 6 & 3+6j & 4+6j \\ 
1 & 3+6j & 4+j & 6 & 3+j & 3+j & 6 & 4+j & 3+6j & 1 & 4+6j & 4+6j  
\end{smallmatrix}
\right] . 
$$

Defining a 1st order pre-addition layer,
we obtain:
\begin{align*}
S_0(1)&=(v_7-v_1),\,S_1(1)=(v_7+v_1),\\ 
S_2(1)&=(v_8-v_2),\,S_3(1)=(v_8+v_2),\\ 
S_4(1)&=(v_9-v_3),\,S_5(1)=(v_9+v_3),\\ 
S_6(1)&=(v_{10}-v_4),\,S_7(1)=(v_{10}+v_4),\\ 
S_8(1)&=(v_{11}-v_5),\,S_9(1)=(v_{11}+v_5),\\ 
S_{10}(1)&=(v_0-v_6),\,S_{11}(1)=(v_0+v_6). 
\end{align*}
Therefore, $V=T^{(1)}S(1)$, where 
$$
T^{(1)}=\left[ {\ 
\begin{smallmatrix}
0 & 1 & 0 & 1 & 0 & 1 & 0 & 1 & 0 & 1 & 0 & 1 \\ 
3+j & 0 & 3+j & 0 & 6 & 0 & 4+j & 0 & 3+6j & 0 & 1 & 0 \\ 
0 & 4+6j & 0 & 3+6j & 0 & 6 & 0 & 3+j & 0 & 4+j & 0 & 1 \\ 
6 & 0 & 1 & 0 & 1 & 0 & 6 & 0 & 6 & 0 & 1 & 0 \\ 
0 & 3+6j & 0 & 3+j & 0 & 1 & 0 & 3+6j & 0 & 3+j & 0 & 1 \\ 
3+6j & 0 & 3+6j & 0 & 6 & 0 & 4+6j & 0 & 3+j & 0 & 1 & 0 \\ 
0 & 6 & 0 & 1 & 0 & 6 & 0 & 1 & 0 & 6 & 0 & 1 \\ 
4+6j & 0 & 3+j & 0 & 1 & 0 & 4+j & 0 & 4+j & 0 & 1 & 0 \\ 
0 & 3+j & 0 & 3+6j & 0 & 1 & 0 & 3+j & 0 & 3+6j & 0 & 1 \\ 
1 & 0 & 1 & 0 & 6 & 0 & 6 & 0 & 1 & 0 & 1 & 0 \\ 
0 & 4+j & 0 & 3+j & 0 & 6 & 0 & 3+6j & 0 & 4+6j & 0 & 1 \\ 
4+j & 0 & 3+6j & 0 & 1 & 0 & 4+6j & 0 & 4+6j & 0 & 1 & 0    
\end{smallmatrix}
}\right] . 
$$

A second order pre-addition layer can be defined according to, 
\begin{align*}
S_0(2)&=S_6(1)-S_0(1),\,S_6(2)=S_9(1)-S_3(1),\\ 
S_1(2)&=S_6(1)+S_0(1),\,S_7(2)=S_9(1)+S_3(1),\\ 
S_2(2)&=S_7(1)-S_1(1),\,S_8(2)=S_{10}(1)-S_4(1),\\ 
S_3(2)&=S_7(1)+S_1(1),\,S_9(2)=S_{10}(1)+S_4(1),\\ 
S_4(2)&=S_8(1)-S_2(1),\,S_{10}(2)=S_{11}(1)-S_5(1),\\ 
S_5(2)&=S_8(1)+S_2(1),\,S_{11}(2)=S_{11}(1)+S_5(1). 
\end{align*}
Therefore, $V=T^{(2)}S(2)$, where 
$$
T^{(2)}=\left[ {\ 
\begin{smallmatrix}
0 & 0 & 0 & 1 & 0 & 0 & 0 & 1 & 0 & 0 & 0 & 1 \\ 
4 & j & 0 & 0 & 6j & 3 & 0 & 0 & 1 & 0 & 0 & 0 \\ 
0 & 0 & 3+j & 0 & 0 & 0 & 4+j & 0 & 0 & 0 & 1 & 0 \\ 
0 & 6 & 0 & 0 & 6 & 0 & 0 & 0 & 0 & 1 & 0 & 0 \\ 
0 & 0 & 0 & 3+6j & 0 & 0 & 0 & 3+j & 0 & 0 & 0 & 1 \\ 
4 & 6j & 0 & 0 & j & 3 & 0 & 0 & 1 & 0 & 0 & 0 \\ 
0 & 0 & 1 & 0 & 0 & 0 & 6 & 0 & 0 & 0 & 1 & 0 \\ 
j & 4 & 0 & 0 & 4 & j & 0 & 0 & 0 & 1 & 0 & 0 \\ 
0 & 0 & 0 & 3+j & 0 & 0 & 0 & 3+6j & 0 & 0 & 0 & 1 \\ 
6 & 0 & 0 & 0 & 0 & 1 & 0 & 0 & 1 & 0 & 0 & 0 \\ 
0 & 0 & 3+6j & 0 & 0 & 0 & 4+6j & 0 & 0 & 0 & 1 & 0 \\ 
j & 4 & 0 & 0 & 4 & 6j & 0 & 0 & 0 & 1 & 0 & 0  
\end{smallmatrix}
}\right] . 
$$
Going further, a 3rd order pre-addition layer is defined: 
\begin{align*}
S_0(3)=&S_4(2)-S_1(2),\,S_6(3)=S_7(2)-S_3(2),\\ 
S_1(3)=&S_4(2)+S_1(2),\,S_7(3)=S_7(2)+S_3(2), \\
S_2(3)=&S_5(2)-S_0(2),\,S_8(3)=S_8(2),\\ 
S_3(3)=&S_5(2)+S_0(2),\,S_9(3)=S_9(2),\\ 
S_4(3)=&S_6(2)-S_2(2),\,S_{10}(3)=S_{10}(2),\\ 
S_5(3)=&S_6(2)+S_2(2),\,S_{11}(3)=S_{11}(2). 
\end{align*}
Thus $V=T^{(3)}S(3)$, where 
$$
T^{(3)}=\left[ 
\begin{smallmatrix}
0 & 0 & 0 & 0 & 0 & 0 & 0 & 1 & 0 & 0 & 0 & 1 \\ 
6j & 0 & 3 & 0 & 0 & 0 & 0 & 0 & 1 & 0 & 0 & 0 \\ 
0 & 0 & 0 & 0 & 4 & j & 0 & 0 & 0 & 0 & 1 & 0 \\ 
0 & 6 & 0 & 0 & 0 & 0 & 0 & 0 & 0 & 1 & 0 & 0 \\ 
0 & 0 & 0 & 0 & 0 & 0 & j & 3 & 0 & 0 & 0 & 1 \\ 
j & 0 & 3 & 0 & 0 & 0 & 0 & 0 & 1 & 0 & 0 & 0 \\ 
0 & 0 & 0 & 0 & 6 & 0 & 0 & 0 & 0 & 0 & 1 & 0 \\ 
0 & 4 & 0 & j & 0 & 0 & 0 & 0 & 0 & 1 & 0 & 0 \\ 
0 & 0 & 0 & 0 & 0 & 0 & 6j & 3 & 0 & 0 & 0 & 1 \\ 
0 & 0 & 1 & 0 & 0 & 0 & 0 & 0 & 1 & 0 & 0 & 0 \\ 
0 & 0 & 0 & 0 & 4 & 6j & 0 & 0 & 0 & 0 & 1 & 0 \\ 
0 & 4 & 6j & 0 & 0 & 0 & 0 & 0 & 0 & 1 & 0 & 0  
\end{smallmatrix}
\right] . 
$$
Since there is just one multiplication by the same factor in columns 2, 3, 5
and 8, the total number of multiplications is 4. The number of additions
required to compute the FFHT is 44 (32 pre-additions and 12 post-additions).
The multiplicative complexity reaches the minimum theoretical complexity and
again the additive complexity is the same as the one obtained for the DHT 
\cite{Renato}.

\subsection{Computing the 16-point FFHT}

Let $v\longleftrightarrow V$ a FFHT\ pair over $GI(p)$. 
Assuming $p=7$ and $\zeta =2+4j$, the corresponding transform is $V=Tv$, where 
$$
T=\left[ 
\begin{smallmatrix}
1 & 1 & 1 & 1 & 1 & 1 & 1 & 1 & 1 & 1 & 1 & 1 & 1 & 1 & 1 & 1 \\ 
1 & 2j & 4 & 2j & 1 & 6j & 0 & j & 6 & 5j & 3 & 2j & 6 & j & 0 & 6j \\ 
1 & 4 & 1 & 0 & 6 & 3 & 6 & 0 & 1 & 4 & 1 & 0 & 6 & 3 & 6 & 0 \\ 
1 & 2j & 0 & 5j & 6 & 6j & 4 & 6j & 6 & 5j & 0 & 2j & 1 & j & 3 & j \\ 
1 & 1 & 6 & 6 & 1 & 1 & 6 & 6 & 1 & 1 & 6 & 6 & 1 & 1 & 6 & 6 \\ 
1 & 6j & 3 & 6j & 1 & 5j & 0 & 2j & 6 & j & 4 & j & 6 & 2j & 0 & 5j \\ 
1 & 0 & 6 & 4 & 6 & 0 & 1 & 3 & 1 & 0 & 6 & 4 & 6 & 0 & 1 & 3 \\ 
1 & 6j & 0 & 6j & 6 & 2j & 3 & 2j & 6 & 6j & 0 & j & 1 & 5j & 4 & 5j \\ 
1 & 6 & 1 & 6 & 1 & 6 & 1 & 6 & 1 & 6 & 1 & 6 & 1 & 6 & 1 & 6 \\ 
1 & 5j & 4 & 5j & 1 & j & 0 & 6j & 6 & 2j & 3 & 2j & 6 & 6j & 0 & j \\ 
1 & 3 & 1 & 0 & 6 & 4 & 6 & 0 & 1 & 3 & 1 & 0 & 6 & 4 & 6 & 0 \\ 
1 & 5j & 0 & 2j & 6 & j & 4 & j & 6 & 2j & 0 & 5j & 1 & 6j & 3 & 6j \\ 
1 & 6 & 6 & 1 & 1 & 6 & 6 & 1 & 1 & 6 & 6 & 1 & 1 & 6 & 6 & 1 \\ 
1 & j & 3 & j & 1 & 2j & 0 & 5j & 6 & 6j & 4 & 6j & 6 & 5j & 0 & 2j \\ 
1 & 0 & 6 & 3 & 6 & 0 & 1 & 4 & 1 & 0 & 6 & 3 & 6 & 0 & 1 & 4 \\ 
1 & 6j & 0 & j & 6 & 5j & 3 & 5j & 6 & j & 0 & 6j & 1 & 2j & 4 & 2j  
\end{smallmatrix}
\right] . 
$$
Now we consider the first order pre-addition layer according to:
\begin{equation}
\label{ord1}
\begin{split}
S_0(1)&=(v_9-v_1),\,S_1(1)=(v_9+v_1),\\ 
S_2(1)&=(v_{10}+v_2),\,S_3(1)=(v_{10}+v_2),\\ 
S_4(1)&=(v_{11}-v_3),\,S_5(1)=(v_{11}+v_3),\\
S_6(1)&=(v_{12}-v_4),\,S_7(1)=(v_{12}+v_4),\\ 
S_8(1)&=(v_{13}-v_5),\,S_9(1)=(v_{13}+v_5),\\ 
S_{10}(1)&=(v_{14}-v_6),\,S_{11}(1)=(v_{14}+v_6),\\ 
S_{12}(1)&=(v_{15}-v_7),\,S_{13}(1)=(v_{15}+v_7),\\ 
S_{14}(1)&=(v_0-v_8),\,S_{15}(1)=(v_0+v_8).
\end{split}
\end{equation}
Therefore, $V=T^{(1)}S(1)$, where 
$$
T^{(1)}=\left[ 
\begin{smallmatrix}
0 & 1 & 0 & 1 & 0 & 1 & 0 & 1 & 0 & 1 & 0 & 1 & 0 & 1 & 0 & 1 \\
5j & 0 & 3 & 0 & 5j & 0 & 6 & 0 & j & 0 & 0 & 0 & 6j & 0 & 1 & 0 \\ 
0 & 4 & 0 & 1 & 0 & 0 & 0 & 6 & 0 & 3 & 0 & 6 & 0 & 0 & 0 & 1 \\ 
5j & 0 & 0 & 0 & 2j & 0 & 1 & 0 & j & 0 & 3 & 0 & j & 0 & 1 & 0 \\ 
0 & 1 & 0 & 6 & 0 & 6 & 0 & 1 & 0 & 1 & 0 & 6 & 0 & 6 & 0 & 1 \\ 
j & 0 & 4 & 0 & j & 0 & 6 & 0 & 2j & 0 & 0 & 0 & 5j & 0 & 1 & 0 \\ 
0 & 0 & 0 & 6 & 0 & 4 & 0 & 6 & 0 & 0 & 0 & 1 & 0 & 3 & 0 & 1 \\ 
6j & 0 & 0 & 0 & j & 0 & 1 & 0 & 5j & 0 & 4 & 0 & 5j & 0 & 1 & 0 \\ 
0 & 6 & 0 & 1 & 0 & 6 & 0 & 1 & 0 & 6 & 0 & 1 & 0 & 6 & 0 & 1 \\ 
2j & 0 & 3 & 0 & 2j & 0 & 6 & 0 & 6j & 0 & 0 & 0 & j & 0 & 1 & 0 \\ 
0 & 3 & 0 & 1 & 0 & 0 & 0 & 6 & 0 & 4 & 0 & 6 & 0 & 0 & 0 & 1 \\ 
2j & 0 & 0 & 0 & 5j & 0 & 1 & 0 & 6j & 0 & 3 & 0 & 6j & 0 & 1 & 0 \\ 
0 & 6 & 0 & 6 & 0 & 1 & 0 & 1 & 0 & 6 & 0 & 6 & 0 & 1 & 0 & 1 \\ 
6j & 0 & 4 & 0 & 6j & 0 & 6 & 0 & 5j & 0 & 0 & 0 & 2j & 0 & 1 & 0 \\ 
0 & 0 & 0 & 6 & 0 & 3 & 0 & 6 & 0 & 0 & 0 & 1 & 0 & 4 & 0 & 1 \\ 
j & 0 & 0 & 0 & 6j & 0 & 1 & 0 & 2j & 0 & 4 & 0 & 2j & 0 & 1 & 0 
\end{smallmatrix}
\right] . 
$$
By defining a 2nd order pre-addition layer, 
we have:
\begin{equation}
\label{ord2}
\begin{split}
S_0(2)&=S_4(1)-S_0(1),\,S_1(2)=S_4(1)+S_0(1),\\
S_2(2)&=S_9(1)-S_1(1),\,S_3(2)=S_9(1)+S_1(1),\\
S_4(2)&=S_2(1),\,S_5(2)=S_{10}(1),\\ 
S_6(2)&=S_{11}(1)-S_3(1),\,S_7(2)=S_{11}(1)+S_3(1),\\ 
S_8(2)&=S_{12}(1)-S_8(1),\,S_9(2)=S_{12}(1)+S_8(1),\\ 
S_{10}(2)&=S_{13}(1)-S_5(1),\,S_{11}(2)=S_{13}(1)+S_5(1),\\ 
S_{12}(2)&=S_{14}(1)-S_6(1),\,S_{13}(2)=S_{14}(1)+S_6(1),\\ 
S_{14}(2)&=S_{15}(1)-S_7(1),\,S_{15}(2)=S_{15}(1)+S_7(1).
\end{split}
\end{equation}
Then, $V=T^{(2)}S(2)$, where 
$$
T^{(2)}=\left[ 
\begin{smallmatrix}
0 & 0 & 0 & 1 & 0 & 0 & 0 & 1 & 0 & 0 & 0 & 1 & 0 & 0 & 0 & 1 \\ 
0 & 5j & 0 & 0 & 3 & 0 & 0 & 0 & 6j & 0 & 0 & 0 & 1 & 0 & 0 & 0 \\ 
0 & 0 & 3 & 0 & 0 & 0 & 6 & 0 & 0 & 0 & 0 & 0 & 0 & 0 & 1 & 0 \\ 
2j & 0 & 0 & 0 & 0 & 3 & 0 & 0 & 0 & j & 0 & 0 & 0 & 1 & 0 & 0 \\ 
0 & 0 & 0 & 1 & 0 & 0 & 0 & 6 & 0 & 0 & 0 & 6 & 0 & 0 & 0 & 1 \\ 
0 & j & 0 & 0 & 4 & 0 & 0 & 0 & 5j & 0 & 0 & 0 & 1 & 0 & 0 & 0 \\ 
0 & 0 & 0 & 0 & 0 & 0 & 1 & 0 & 0 & 0 & 3 & 0 & 0 & 0 & 1 & 0 \\ 
j & 0 & 0 & 0 & 0 & 4 & 0 & 0 & 0 & 5j & 0 & 0 & 0 & 1 & 0 & 0 \\ 
0 & 0 & 0 & 6 & 0 & 0 & 0 & 1 & 0 & 0 & 0 & 6 & 0 & 0 & 0 & 1 \\ 
0 & 2j & 0 & 0 & 3 & 0 & 0 & 0 & j & 0 & 0 & 0 & 1 & 0 & 0 & 0 \\ 
0 & 0 & 4 & 0 & 0 & 0 & 6 & 0 & 0 & 0 & 0 & 0 & 0 & 0 & 1 & 0 \\ 
5j & 0 & 0 & 0 & 0 & 3 & 0 & 0 & 0 & 6j & 0 & 0 & 0 & 1 & 0 & 0 \\ 
0 & 0 & 0 & 6 & 0 & 0 & 0 & 6 & 0 & 0 & 0 & 1 & 0 & 0 & 0 & 1 \\ 
0 & 6j & 0 & 0 & 4 & 0 & 0 & 0 & 2j & 0 & 0 & 0 & 1 & 0 & 0 & 0 \\ 
0 & 0 & 0 & 0 & 0 & 0 & 1 & 0 & 0 & 0 & 4 & 0 & 0 & 0 & 1 & 0 \\ 
6j & 0 & 0 & 0 & 0 & 4 & 0 & 0 & 0 & 2j & 0 & 0 & 0 & 1 & 0 & 0  
\end{smallmatrix}
\right] . 
$$
The columns do not cope. However, multiplying both the 5 and 6 columns by 
$2\,\in GF(7)$, they can be combined with columns 13 and 14, respectively.
Defining then a 2nd order pre-addition layer (with two multiplications in
columns 10 and 11), 
we have:
\begin{align*}
S_0(3)&=S_0(2),\,S_1(3)=S_1(2),\\
S_2(3)&=S_8(2),\,S_3(3)=S_9(2),\\
S_4(3)&=S_2(2),\,S_5(3)=S_{10}(2),\\
S_6(3)&=S_{11}(2)-S_3(2),\,S_7(3)=S_{11}(2)+S_3(2),\\ 
S_8(3)&=S_{12}(2)-4S_4(2),\,S_9(3)=S_{12}(2)+2S_4(2),\\ 
S_{10}(3)&=S_{13}(2)-4S_5(2),S_{11}(3)=S_{13}(2)+2S_5(2),\\
S_{12}(3)&=S_{14}(2)-S_6(2),\,S_{13}(3)=S_{14}(2)+S_6(2),\\ 
S_{14}(3)&=S_{15}(2)-S_7(2),\,S_{15}(3)=S_{15}(2)+S_7(2).
\end{align*}
Finally, $V=T^{(3)}S(3)$, where 
$$
T^{(3)}=\left[ 
\begin{smallmatrix}
0 & 0 & 0 & 0 & 0 & 0 & 0 & 1 & 0 & 0 & 0 & 0 & 0 & 0 & 0 & 1 \\ 
0 & 5j & 6j & 0 & 0 & 0 & 0 & 0 & 1 & 0 & 0 & 0 & 0 & 0 & 0 & 0 \\ 
0 & 0 & 0 & 0 & 3 & 0 & 0 & 0 & 0 & 0 & 0 & 0 & 1 & 0 & 0 & 0 \\ 
2j & 0 & 0 & j & 0 & 0 & 0 & 0 & 0 & 0 & 1 & 0 & 0 & 0 & 0 & 0 \\ 
0 & 0 & 0 & 0 & 0 & 0 & 6 & 0 & 0 & 0 & 0 & 0 & 0 & 0 & 1 & 0 \\ 
0 & j & 5j & 0 & 0 & 0 & 0 & 0 & 0 & 1 & 0 & 0 & 0 & 0 & 0 & 0 \\ 
0 & 0 & 0 & 0 & 0 & 3 & 0 & 0 & 0 & 0 & 0 & 0 & 0 & 1 & 0 & 0 \\ 
j & 0 & 0 & 5j & 0 & 0 & 0 & 0 & 0 & 0 & 0 & 1 & 0 & 0 & 0 & 0 \\ 
0 & 0 & 0 & 0 & 0 & 0 & 0 & 6 & 0 & 0 & 0 & 0 & 0 & 0 & 0 & 1 \\ 
0 & 2j & j & 0 & 0 & 0 & 0 & 0 & 1 & 0 & 0 & 0 & 0 & 0 & 0 & 0 \\ 
0 & 0 & 0 & 0 & 4 & 0 & 0 & 0 & 0 & 0 & 0 & 0 & 1 & 0 & 0 & 0 \\ 
5j & 0 & 0 & 6j & 0 & 0 & 0 & 0 & 0 & 0 & 1 & 0 & 0 & 0 & 0 & 0 \\ 
0 & 0 & 0 & 0 & 0 & 0 & 1 & 0 & 0 & 0 & 0 & 0 & 0 & 0 & 1 & 0 \\ 
0 & 6j & 2j & 0 & 0 & 0 & 0 & 0 & 0 & 1 & 0 & 0 & 0 & 0 & 0 & 0 \\ 
0 & 0 & 0 & 0 & 0 & 4 & 0 & 0 & 0 & 0 & 0 & 0 & 0 & 1 & 0 & 0 \\ 
6j & 0 & 0 & 2j & 0 & 0 & 0 & 0 & 0 & 0 & 0 & 1 & 0 & 0 & 0 & 0 
\end{smallmatrix}
\right] . 
$$
There is only one multiplication in columns 1, 2, 3, 4, 5 and 6, besides
two pre-multiplications, so the total complexity is 8. The number of
additions is 50 (40 pre-additions and 16 post-additions).
In this case, the
number of multiplications is less than 10, the minimum expected
multiplication complexity~\cite{Complexid}. It can be concluded that there
are two trivial multiplications. It is not simple to identify which are the
trivial multiplications from the observation of matrices over $GI(7)$.
Carrying on the same analysis over another finite field, the same
combination of columns was observed, i.e. the approach does not depend on
the finite field but on the length. Let $v\longleftrightarrow V$ be an FFHT\
pair over $GI(p)$. Considering now $p=31,\zeta =7+13j$, the transform matrix 
$T$ will be, 
$$
T=\left[ 
\begin{smallmatrix}
1 & 1 & 1 & 1 & 1 & 1 & 1 & 1 & 1 & 1 & 1 & 1 & 1 & 1 & 1 & 1 \\ 
1 & 20 & 0 & 11 & 30 & 6 & 23 & 6 & 30 & 11 & 0 & 20 & 1 & 25 & 8 & 25 \\ 
1 & 0 & 30 & 23 & 30 & 0 & 1 & 8 & 1 & 0 & 30 & 23 & 30 & 0 & 1 & 8 \\ 
1 & 11 & 23 & 11 & 1 & 25 & 0 & 6 & 30 & 20 & 8 & 20 & 30 & 6 & 0 & 25 \\ 
1 & 30 & 30 & 1 & 1 & 30 & 30 & 1 & 1 & 30 & 30 & 1 & 1 & 30 & 30 & 1 \\ 
1 & 6 & 0 & 25 & 30 & 11 & 8 & 11 & 30 & 25 & 0 & 6 & 1 & 20 & 23 & 20 \\ 
1 & 23 & 1 & 0 & 30 & 8 & 30 & 0 & 1 & 23 & 1 & 0 & 30 & 8 & 30 & 0 \\ 
1 & 6 & 8 & 6 & 1 & 11 & 0 & 20 & 30 & 25 & 23 & 25 & 30 & 20 & 0 & 11 \\ 
1 & 30 & 1 & 30 & 1 & 30 & 1 & 30 & 1 & 30 & 1 & 30 & 1 & 30 & 1 & 30 \\ 
1 & 11 & 0 & 20 & 30 & 25 & 23 & 25 & 30 & 20 & 0 & 11 & 1 & 6 & 8 & 6 \\ 
1 & 0 & 20 & 8 & 30 & 0 & 1 & 23 & 1 & 0 & 30 & 8 & 30 & 0 & 1 & 23 \\ 
1 & 20 & 23 & 20 & 1 & 6 & 0 & 25 & 30 & 11 & 8 & 11 & 30 & 25 & 0 & 6 \\ 
1 & 1 & 30 & 30 & 1 & 1 & 30 & 30 & 1 & 1 & 30 & 30 & 1 & 1 & 30 & 30 \\ 
1 & 25 & 0 & 6 & 30 & 20 & 8 & 20 & 30 & 6 & 0 & 25 & 1 & 11 & 23 & 11 \\ 
1 & 8 & 1 & 0 & 30 & 23 & 30 & 0 & 1 & 8 & 1 & 0 & 30 & 23 & 30 & 0 \\ 
1 & 25 & 8 & 25 & 1 & 20 & 0 & 11 & 30 & 6 & 23 & 6 & 30 & 11 & 0 & 20   
\end{smallmatrix}
\right] . 
$$
Considering the first order pre-addition layer (Equation~\ref{ord1}) and the
second order pre-addition layer (Equation~\ref{ord2}), $V=T^{(2)}S(2)$, where 
$$
T^{(2)}=\left[ 
\begin{smallmatrix}
0 & 0 & 0 & 1 & 0 & 0 & 0 & 1 & 0 & 0 & 0 & 1 & 0 & 0 & 0 & 1 \\ 
20 & 0 & 0 & 0 & 0 & 8 & 0 & 0 & 0 & 25 & 0 & 0 & 0 & 1 & 0 & 0 \\ 
0 & 0 & 0 & 0 & 0 & 0 & 1 & 0 & 0 & 0 & 8 & 0 & 0 & 0 & 1 & 0 \\ 
0 & 20 & 0 & 0 & 8 & 0 & 0 & 0 & 25 & 0 & 0 & 0 & 1 & 0 & 0 & 0 \\ 
0 & 0 & 0 & 30 & 0 & 0 & 0 & 30 & 0 & 0 & 0 & 1 & 0 & 0 & 0 & 1 \\ 
6 & 0 & 0 & 0 & 0 & 23 & 0 & 0 & 0 & 20 & 0 & 0 & 0 & 1 & 0 & 0 \\ 
0 & 0 & 8 & 0 & 0 & 0 & 30 & 0 & 0 & 0 & 0 & 0 & 0 & 0 & 1 & 0 \\ 
0 & 25 & 0 & 0 & 23 & 0 & 0 & 0 & 11 & 0 & 0 & 0 & 1 & 0 & 0 & 0 \\ 
0 & 0 & 0 & 30 & 0 & 0 & 0 & 1 & 0 & 0 & 0 & 30 & 0 & 0 & 0 & 1 \\ 
11 & 0 & 0 & 0 & 0 & 8 & 0 & 0 & 0 & 6 & 0 & 0 & 0 & 1 & 0 & 0 \\ 
0 & 0 & 0 & 0 & 0 & 0 & 1 & 0 & 0 & 0 & 23 & 0 & 0 & 0 & 1 & 0 \\ 
0 & 11 & 0 & 0 & 8 & 0 & 0 & 0 & 6 & 0 & 0 & 0 & 1 & 0 & 0 & 0 \\ 
0 & 0 & 0 & 1 & 0 & 0 & 0 & 30 & 0 & 0 & 0 & 30 & 0 & 0 & 0 & 1 \\ 
25 & 0 & 0 & 0 & 0 & 23 & 0 & 0 & 0 & 11 & 0 & 0 & 0 & 1 & 0 & 0 \\ 
0 & 0 & 23 & 0 & 0 & 0 & 30 & 0 & 0 & 0 & 0 & 0 & 0 & 0 & 1 & 0 \\ 
0 & 6 & 0 & 0 & 23 & 0 & 0 & 0 & 20 & 0 & 0 & 0 & 1 & 0 & 0 & 0 
\end{smallmatrix}
\right] . 
$$
Again, some columns do not cope. But, multiplying both columns 5 and 6 by 
$4\,\in GF(31)$, they can be combined with columns 14 and 13, respectively. The
same occurs with columns 9 and 10, which 
combine with columns 2 and 1,
respectively.
A third layer of pre-additions, including four
pre-multiplications in columns 1, 2, 5, and 6, 
is given by:
\begin{align*}
S_0(3)&=S_7(2)-S_3(2),\,S_1(3)=S_7(2)+S_7(2),\\
S_2(3)&=S_8(2)-7S_1(2),\, S_3(3)=S_7(2)+8S_1(2),\\
S_4(3)&=S_9(2)-7S_0(2),\,S_5(3)=S_9(2)+7S_0(2),\\
S_6(3)&=S_2(2),\,S_7(3)=S_{10}(2),\\ 
S_8(3)&=S_{12}(2)-4S_4(2),\,S_9(3)=S_{12}(2)+4S_4(2),\\ 
S_{10}(3)&=S_{13}(2)-4S_5(2),\,S_{11}(3)=S_{13}(2)+4S_5(2),\\
S_{12}(3)&=S_{14}(2)-S_6(2),\,S_{13}(3)=S_{14}(2)+S_6(2),\\
S_{14}(3)&=S_{15}(2)-S_{11}(2),\,S_{15}(3)=S_{15}(2)+S_{11}(2).
\end{align*}
Therefore, 
we have that $V=T^{(3)}S(3)$, where 
$$
T^{(3)}=\left[ 
\begin{smallmatrix}
0 & 1 & 0 & 0 & 0 & 0 & 0 & 0 & 0 & 0 & 0 & 0 & 0 & 0 & 0 & 1 \\ 
0 & 0 & 0 & 0 & 0 & 20 & 0 & 0 & 0 & 0 & 0 & 1 & 0 & 0 & 0 & 0 \\ 
0 & 0 & 0 & 0 & 0 & 0 & 0 & 8 & 0 & 0 & 0 & 0 & 0 & 1 & 0 & 0 \\ 
0 & 0 & 0 & 20 & 0 & 0 & 0 & 0 & 0 & 1 & 0 & 0 & 0 & 0 & 0 & 0 \\ 
0 & 30 & 0 & 0 & 0 & 0 & 0 & 0 & 0 & 0 & 0 & 0 & 0 & 0 & 0 & 1 \\ 
0 & 0 & 0 & 0 & 0 & 6 & 0 & 0 & 0 & 0 & 1 & 0 & 0 & 0 & 0 & 0 \\ 
0 & 0 & 0 & 0 & 0 & 0 & 8 & 0 & 0 & 0 & 0 & 0 & 1 & 0 & 0 & 0 \\ 
0 & 0 & 0 & 25 & 0 & 0 & 0 & 0 & 1 & 0 & 0 & 0 & 0 & 0 & 0 & 0 \\ 
1 & 0 & 0 & 0 & 0 & 0 & 0 & 0 & 0 & 0 & 0 & 0 & 0 & 0 & 1 & 0 \\ 
0 & 0 & 0 & 0 & 0 & 11 & 0 & 0 & 0 & 0 & 0 & 1 & 0 & 0 & 0 & 0 \\ 
0 & 0 & 0 & 0 & 0 & 0 & 0 & 23 & 0 & 0 & 0 & 0 & 0 & 1 & 0 & 0 \\ 
0 & 0 & 0 & 11 & 0 & 0 & 0 & 0 & 0 & 1 & 0 & 0 & 0 & 0 & 0 & 0 \\ 
30 & 0 & 0 & 0 & 0 & 0 & 0 & 0 & 0 & 0 & 0 & 0 & 0 & 0 & 1 & 0 \\ 
0 & 0 & 0 & 0 & 0 & 25 & 0 & 0 & 0 & 0 & 1 & 0 & 0 & 0 & 0 & 0 \\ 
0 & 0 & 0 & 0 & 0 & 0 & 23 & 0 & 0 & 0 & 0 & 0 & 1 & 0 & 0 & 0 \\ 
0 & 0 & 0 & 6 & 0 & 0 & 0 & 0 & 1 & 0 & 0 & 0 & 0 & 0 & 0 & 0
\end{smallmatrix}
\right] . 
$$
The required number of multiplications is 10 (4 pre-multiplications, 2
multiplications in column 4, 2 multiplications in column 5, 1 multiplication in column 7
and 1 multiplication in column 8). The number of
additions is 60 (44 pre-additions, 16 post-additions). A complexity
comparison of $N$-point FFHT fast algorithms (for $N=8$ and $N=16$) 
is given in
Tables~\ref{tabela1} and~\ref{tabela2}.

\begin{table}
\begin{center}
\caption{Complexity of the $8$-point FFHT}
\begin{tabular}{ccccccc}
\toprule
Fast algorithms & $M(8)$ & $A(8)$ & $ M(8)+A(8)$  \\
\midrule
Cooley-Tukey-4  & 12 & 48 & 60 \\
Split-Radix & 8 & 42 & 50  \\
Cooley-Tukey-2 & 4 & 26 & 30 \\
Rader-Brenner & 2  & 24 & 26  \\
Proposed & 2 & 22 & 24 \\
\bottomrule
\end{tabular}
\label{tabela1}

\end{center}
\end{table}
\begin{table}
\caption{Complexity of the $16$-point FFHT}
\begin{center}
\begin{tabular}{ccccccc}
\toprule
Fast algorithms &  $M(16)$ & $A(16)$ & 
$M(16)+A(16)$ \\
\midrule
Cooley-Tukey-2 & 20 & 74 & 94 \\
Cooley-Tukey-4  & 14 & 70  & 84 \\
Split-Radix &  12 & 64 & 76 \\
Rader-Brenner & 10 & 64 & 74 \\
Proposed & 10 & 60 & 70\\
\bottomrule
\end{tabular}
\label{tabela2}
\end{center}
\end{table}

In the above examples, the Hadamard
decomposition algorithm presents a lower complexity to compute an FFHT
compared to existing FFFT/DFT\ algorithms. Multiplicative complexity saving
regarding classical Cooley-Tukey is 50\% ($N=16$ and $N=8$). The total complexity
saving regarding the same algorithm is roughly \linebreak 25\% ($N=16$), 20\% ($N=8$).

\section{Conclusions}

Fast algorithms
for the finite field Hartley transform 
based on Walsh-Hadamard decompositions
were developed and applied to short blocklengths.
The theoretical multiplicative complexity lower bounds
were achieved.
The total complexity (additive and multiplicative)
of the algorithms was compared to
that of popular algorithms and
the lower values were obtained for the FT.
These FTs are  
attractive
and easy to implement using low-cost high-speed dedicated 
circuitry.

\section*{Acknowledgments}

This work was partially supported by CNPq and CAPES.

\end{document}